\documentclass[reqno, 11pt, a4paper]{amsart}

\usepackage[utf8]{inputenc}
\usepackage[T1]{fontenc}
\usepackage{lmodern}
\usepackage{amssymb}
\usepackage[colorlinks=true,citecolor=blue,urlcolor=blue,linkcolor=blue,ocgcolorlinks]{hyperref}
\usepackage{graphicx}
\usepackage{tikz-cd}
\usepackage{mathrsfs}
\usepackage{bbm}
\usepackage{ytableau}

\makeatletter
\renewcommand*{\@secnumfont}{\bfseries}
\def\section{\@startsection{section}{1}%
  \z@{.7\linespacing\@plus\linespacing}{.5\linespacing}%
  {\normalfont\bfseries\centering}}
\def\subsection{\@startsection{subsection}{2}%
  \z@{.5\linespacing\@plus.7\linespacing}{-.5em}%
  {\normalfont\itshape}}
\makeatother

\def\secskip{\vspace{.5\linespacing plus.7\linespacing}}

\newcommand*{\bC}{\mathbb{C}}

\newcommand*{\bL}{\mathbb{L}}

\newcommand*{\bQ}{\mathbb{Q}}

\newcommand*{\sM}{\mathscr{M}}

\newcommand*{\fS}{\mathfrak{S}}
\newcommand*{\ssS}{\mathsf{S}}
\newcommand*{\bone}{\mathbbm{1}}

\newcommand*{\eg}{\textit{e.g.}~}

\newcommand*{\ie}{\textit{i.e.}~}
\newcommand*{\loccit}{\textit{loc.\,cit.}}

\newcommand*{\rand}{\textup{ \ and \ }}

\newcommand*{\rin}{\textup{ \ in \ }}

\newcommand*{\rleft}{\textup{(}}
\newcommand*{\rright}{\textup{)}}

\newcommand*{\alg}{\mathop\mathrm{alg}\nolimits}

\newcommand*{\cl}{\mathop\mathrm{cl}\nolimits}

\newcommand*{\CH}{\mathop\mathrm{CH}\nolimits}

\newcommand*{\even}{\mathop\mathrm{even}\nolimits}

\newcommand*{\Mon}{\mathop\mathrm{Mon}\nolimits}

\newcommand*{\odd}{\mathop\mathrm{odd}\nolimits}

\newcommand*{\Pic}{\mathop\mathrm{Pic}\nolimits}
\newcommand*{\pr}{\mathop\mathrm{pr}\nolimits}

\newcommand*{\rt}{\mathop\mathrm{rt}\nolimits}
\newcommand*{\sgn}{\mathop\mathrm{sgn}\nolimits}

\newcommand*{\Tab}{\mathop\mathrm{Tab}\nolimits}
\newcommand*{\tr}{\mathop\mathrm{tr}\nolimits}

\newcommand*{\tofrom}{\longleftrightarrow}

\makeatletter
\newcommand*{\vast}{\bBigg@{3}}
\newcommand*{\Vast}{\bBigg@{4}}
\makeatother

\begin{document}

\title{Finite-dimensionality and cycles on powers of $K3$ surfaces}
\author{Qizheng Yin}
\date{\today}
\address{\rm ETH Zürich, Departement Mathematik, Rämistrasse 101, 8092 Zürich, Switzerland}
\email{\href{mailto:yin@math.ethz.ch}{yin@math.ethz.ch}}

\begin{abstract}
For a $K3$ surface $S$, consider the subring of $\CH(S^n)$ generated by divisor and diagonal classes (with $\bQ$-coefficients). Voisin conjectures that the restriction of the cycle class map to this ring is injective. We prove that Voisin's conjecture is equivalent to the finite-dimensionality of $S$ in the sense of Kimura-O'Sullivan. As a consequence, we obtain examples of $S$ whose Hilbert schemes satisfy the Beauville-Voisin conjecture.
\end{abstract}
\maketitle

\vspace{-1cm}
\vfill

\section*{Introduction}
\secskip

\noindent Let $S$ be a smooth projective $K3$ surface over a field $k$. For $n \geq 1$, consider the $\bQ$-subalgebra of the Chow ring $\CH(S^n)$ with $\bQ$-coefficients generated by (pull-backs of) divisor classes on $S$ and the diagonal class on $S \times S$. We denote it by $R(S^n)$. Regarding its structure, Voisin made the following conjecture (\cite{Voi08}, Conjecture~1.6).

\subsection*{Conjecture~1} \label{conj:v}--- \textit{For $n \geq 1$, the restriction of the cycle class map $\cl \colon \CH(S^n) \to H(S^n)$ to $R(S^n)$ is injective.}

\secskip
The case $n = 1$ is the well-known result of Beauville and Voisin (\cite{BV04}, Theorem~1). Voisin also proved Conjecture~\hyperref[conj:v]{1} for $n \leq 2b_{\tr} + 1$, where $b_{\tr} = 22 - \rho$ is the rank of the transcendental part of $H^2(S)$ (\cite{Voi08}, Proposition~2.2). As is remarked in \cite{Voi14}, Section~5.1, Conjecture~\hyperref[conj:v]{1} turns out to be rather strong. Notably it implies the finite-dimensionality of $S$ in the sense of Kimura-O'Sullivan (\cite{Kim05}; see also Section~\ref{sec:kimura}), a conjecture that is widely open even for $K3$ surfaces. 

The aim of this short note is to prove the converse: that finite-dimensionality suffices to deduce Conjecture~\hyperref[conj:v]{1}.

\subsection*{Theorem} --- \textit{Conjecture~\hyperref[conj:v]{1} holds for $S$ if and only if $S$ is finite-dimensional.}

\secskip
More precisely, we prove that the relations found in \cite{BV04} (see also Section~\ref{sec:bv}) plus the one predicted by the finite-dimensionality of $S$ generate all relations in $\cl\big(R(S^n)\big)$. To achieve this, we reduce the problem to a manageable algebraic form, whose solution has long been known to algebraists (Hanlon and Wales \cite{HW89}).

Further, let $k$ be algebraically closed. As is shown by Voisin (\cite{Voi08}, Proposition~2.5), Conjecture~\hyperref[conj:v]{1} implies the following conjecture for the Hilbert schemes $S^{[n]}$ of $S$ (Conjecture~1.3 in \loccit; stated for $k = \bC$ and often referred to as the Beauville-Voisin conjecture).

\subsection*{Conjecture~2} \label{conj:bv} --- \textit{Let $X$ be an irreducible holomorphic symplectic variety \rleft hyper-Kähler manifold\rright. Then the restriction of the cycle class map $\cl \colon \CH(X) \to H(X)$ to the $\bQ$-subalgebra generated by divisor classes and Chern classes of the tangent bundle is injective.} 

\secskip
We thus obtain an immediate consequence.

\subsection*{Corollary} --- \textit{If $S$ is finite-dimensional, then Conjecture~\hyperref[conj:bv]{2} holds for $S^{[n]}$ and for all $n \geq 1$.}

\newpage
Finally, we refer to \cite{Ped13} in characteristic~$0$ and \cite{Lie13} in positive characteristic for $K3$ surfaces known to be finite-dimensional. Among them are Kummer surfaces\footnote{For Kummer surfaces, one can also prove Conjectures~\hyperref[conj:v]{1} and~\hyperref[conj:bv]{2} by applying results on abelian varieties (\eg \cite{Moo11}, Corollary 9.4). We refer to \cite{Fu13} for a similar argument proving Conjecture~\hyperref[conj:bv]{2} for generalized Kummer varieties.}, surfaces of Picard rank $19, 20$ and $22$ (supersingular), and some sporadic cases of even Picard rank.

\subsection*{Notation and conventions} --- Throughout, Chow groups $\CH$ are with $\bQ$-coefficients. We fix a Weil cohomology theory $H$ (\eg singular cohomology when $k = \bC$, or $\ell$-adic cohomology in general).

\subsection*{Acknowledgements} --- Thanks to Mehdi Tavakol for explaining the result of Hanlon and Wales, and to Claire Voisin and Rahul Pandharipande for useful discussions. This work was carried out in the group of Pandharipande at ETH Zürich, supported by grant ERC-2012-AdG-320368-MCSK.

\secskip
\section{Preliminaries}
\secskip

\subsection{\texorpdfstring{\!\!\!}{}} Let $S$ be a smooth projective $K3$ surface over $k$. Denote by $o \in \CH^2(S)$ the class of a distinguished point on $S$, as in \cite{BV04}, Theorem~1. Take a basis $\{L^s\}_{1 \leq s \leq \rho}$ of $\Pic(S)$, and write $l^s = c_1(L^s) \in \CH^1(S)$. For convenience we assume $\{l^s\}$ to be orthogonal. Further, denote by $\delta = [\Delta] \in \CH^2(S \times S)$ the class of the diagonal.

Consider the projections $\pr_i \colon S^n \to S$ for $1 \leq i \leq n$, and $(\pr_i, \pr_j) \colon S^n \to S \times S$ for $1 \leq i, j \leq n$ and $i \neq j$. We write
\begin{equation*}
o_i = \pr_i^*(o) \in \CH^2(S^n), \ \ l^s_i = \pr_i^*(l^s) \in \CH^1(S^n), \rand \delta_{i, j} = (\pr_i, \pr_j)^*(\delta) \in \CH^2(S^n).
\end{equation*}
The ring $R(S^n)$ is defined to be the $\bQ$-subalgebra of $\CH(S^n)$ generated by $\{o_i\}, \{l^s_i\}$ and $\{\delta_{i, j}\}$.

\subsection{Relations in \texorpdfstring{$R(S^n)$}{R(Sn)}} \label{sec:bv} --- The following set of relations summarizes the main results of \cite{BV04}\footnote{By lifting to characteristic~$0$, the results of \cite{BV04} remain valid in positive characteristic.}.
\begin{gather}
o_i \cdot o_i = 0, \ \ l^s_i \cdot o_i = 0, \rand l^s_i \cdot l^s_i = \deg(l^s_i)^2 o_i; \label{eq:bv1} \\
\delta_{i, j} \cdot o_i = o_i \cdot o_j, \ \ \delta_{i, j} \cdot l^s_i = l^s_i \cdot o_j + o_i \cdot l^s_j, \rand \delta_{i, j} \cdot \delta_{i, j} = 24o_i \cdot o_j; \label{eq:bv2} \\
\delta_{i, j} \cdot \delta_{i, k} = \delta_{i, j} \cdot o_k +  \delta_{i, k} \cdot o_j + \delta_{j, k} \cdot o_i - o_i \cdot o_j - o_i \cdot o_k - o_j \cdot o_k. \label{eq:bv3}
\end{gather}
Note that the relations in \eqref{eq:bv1}, \eqref{eq:bv2} and \eqref{eq:bv3} involve $1, 2$ and $3$ factors of $S^n$ respectively.

As we will see, it is both meaningful and convenient to replace $\delta_{i, j}$ (as generator of $R(S^n)$) by
\begin{equation*}
\tau_{i, j} = \delta_{i, j} - o_i - o_j - \sum_{s = 1}^{\rho}\frac{l^s_i \cdot l^s_j}{\deg(l^s)^2} \in \CH^2(S^n).
\end{equation*}
Here $\tau$ stands for ``transcendental''. The relations above now appear in an even simpler form.

\subsection{Lemma} --- {\it In $R(S^n)$ we have relations
\begin{gather}
o_i \cdot o_i = 0, \ \ l^s_i \cdot o_i = 0, \rand l^s_i \cdot l^s_i = \deg(l^s_i)^2 o_i; \label{eq:bv1'} \\
\tau_{i, j} \cdot o_i = 0, \ \ \tau_{i, j} \cdot l^s_i = 0, \rand \tau_{i, j} \cdot \tau_{i, j} = b_{\tr} o_i \cdot o_j; \label{eq:bv2'} \\
\tau_{i, j} \cdot \tau_{i, k} = \tau_{j, k} \cdot o_i, \label{eq:bv3'}
\end{gather}
where $b_{\tr} = 22 - \rho$ is the rank of the transcendental part of $H^2(S)$.}

\begin{proof}
The calculation is straightforward and we only do \eqref{eq:bv3'}. By \eqref{eq:bv2}, \eqref{eq:bv3} and \eqref{eq:bv2'}, we get
\begin{align*}
\tau_{i, j} \cdot \tau_{i, k} & = \tau_{i, j} \cdot \bigg(\delta_{i, k} - o_i - o_k - \sum_{s = 1}^{\rho}\frac{l^s_i \cdot l^s_k}{\deg(l^s)^2}\bigg) \displaybreak \\
& = \tau_{i, j} \cdot \delta_{i, k} - \tau_{i, j} \cdot o_k \\
& = \bigg(\delta_{i, j} - o_i - o_j - \sum_{s = 1}^{\rho}\frac{l^s_i \cdot l^s_j}{\deg(l^s)^2}\bigg) \cdot \delta_{i, k} - \bigg(\delta_{i, j} - o_i - o_j - \sum_{s = 1}^{\rho}\frac{l^s_i \cdot l^s_j}{\deg(l^s)^2}\bigg) \cdot o_k \\
& = \delta_{i, j} \cdot \delta_{i, k} - o_i \cdot o_k - \delta_{i, k} \cdot o_j - \delta_{i, j} \cdot o_k + o_i \cdot o_k + o_j \cdot o_k - \bigg(\sum_{s = 1}^{\rho}\frac{l^s_j \cdot l^s_k}{\deg(l^s)^2}\bigg) \cdot o_i \\
& = \delta_{j, k} \cdot o_i - o_i \cdot o_j - o_i \cdot o_k - \bigg(\sum_{s = 1}^{\rho}\frac{l^s_j \cdot l^s_k}{\deg(l^s)^2}\bigg) \cdot o_i = \tau_{j, k} \cdot o_i. \qedhere
\end{align*}
\end{proof}

\subsection{Finite-dimensionality} \label{sec:kimura} --- We refer to \cite{And04}, Chapitre~4 for the definition of Chow motives over $k$. A motive $M$ is said to be finite-dimensional if $M$ can be decomposed into $M^{\odd} \oplus M^{\even}$ satisfying $\ssS^{N_1}(M^{\odd}) = 0$ and $\wedge^{N_2}(M^{\even}) = 0$ for some $N_1, N_2 > 0$. Here $\ssS$ and $\wedge$ are the symmetric and exterior powers respectively. More precisely, if $M$ is finite-dimensional, one can take $N_1 = \dim H^{\odd}(M) + 1$ and $N_2 = \dim H^{\even}(M) + 1$.

It is conjectured that all Chow motives are finite-dimensional (\cite{Kim05}, Conjecture~7.1), although this is proven only for the subcategory generated by the motives of curves (Theorem~4.2 in \loccit). The motive of a $K3$ surface is believed to be in this category (\eg over $\bC$ by applying the Kuga-Satake construction and the Lefschetz standard conjecture in addition). However, as we discussed, even its finite-dimensionality remains unknown in general.

\subsection{\texorpdfstring{\!\!\!}{}} Back to the $K3$ surface $S$. We now interpret what it means for $S$ to be finite-dimensional. By \cite{KMP07}, Section~7.2, the motive of $S$ (denoted by $h(S)$) admits a decomposition
\begin{equation*}
h(S) = h^0(S) \oplus h^2_{\alg}(S) \oplus h^2_{\tr}(S) \oplus h^4(S) = \bone \oplus \bL^{\oplus \rho} \oplus h^2_{\tr}(S) \oplus \bL^{\otimes2}.
\end{equation*}
Here $\bone$ is the unit motive and $\bL$ is the Lefschetz motive, both of which are (evenly) finite-dimensional. The only part that remains unclear is the motive $h^2_{\tr}(S)$, which is define exactly by the projector $\tau = \tau_{1, 2} \in \CH^2(S \times S)$. We have $\dim H\big(h^2_{\tr}(S)\big) = \dim H^2\big(h^2_{\tr}(S)\big) = b_{\tr}$. It follows that $S$ is finite-dimensional if and only if $\wedge^{b_{\tr} + 1}h^2_{\tr}(S) = 0$. In down-to-earth terms, this means
\begin{equation} \label{eq:motrel}
\sum_{g \in \fS_{b_{\tr} + 1}} \sgn(g) \prod_{i = 1}^{b_{\tr} + 1} \tau_{i, b_{\tr} + 1 + g(i)} = 0 \rin R^{b_{\tr} + 1}\big(S^{2(b_{\tr} + 1)}\big),
\end{equation} 
where $\fS$ stands for the symmetric group and $\sgn$ the signature. Note that \eqref{eq:motrel} holds in $H\big(S^{2(b_{\tr} + 1)}\big)$ by definition, so Conjecture~\hyperref[conj:v]{1} implies the finite-dimensionality of $S$.

The group $\fS_{2(b_{\tr} + 1)}$ acts on $S^{2(b_{\tr} + 1)}$ by permutations. It then acts on \eqref{eq:motrel} and produces more relations in $R\big(S^{2(b_{\tr} + 1)}\big)$.

\section{Proof of the theorem}
\secskip

\subsection{\texorpdfstring{\!\!\!}{}} Using the relations \eqref{eq:bv1'}, \eqref{eq:bv2'} and \eqref{eq:bv3'}, it is not difficult to see that  $R^{2n}(S^n) = \bQ \cdot \{\prod_{i = 1}^n o_i\}$. Then for $0 \leq m \leq 2n$, consider the pairing between $R^m(S^n)$ and $R^{2n - m}(S^n)$. We will show that by assuming \eqref{eq:motrel} and its permutations, the pairing is already perfect. This means there cannot be more relations in $\cl\big(R(S^n)\big)$ than in $R(S^n)$, which proves the theorem.

\subsection{\texorpdfstring{\!\!\!}{}} The first step is basically the same as in \cite{Voi08}, proof of Lemma~2.3. By applying \eqref{eq:bv1'}, \eqref{eq:bv2'} and \eqref{eq:bv3'}, one observes that $R(S^n)$ is linearly spanned by monomials in $\{o_i\}, \{l^s_i\}$ and $\{\tau_{i, j}\}$ with no repeated index, \ie each index $i \in \{1, \ldots, n\}$ appears at most once. We denote by $\Mon^m(n)$ the set of all such monomials that lie in $R^m(S^n)$ (when $m = 0$ we set $\Mon^0(n) = \{1\}$).

An element in $\Mon^m(n)$ can be uniquely written as
\begin{equation*}
\tau_{I, \alpha} \cdot l_{J, \beta} \cdot o_K.
\end{equation*}
Here $I, J$ and $K$ are pairwise disjoint subsets of $\{1, \ldots, n\}$ satisfying $|I| + |J| + 2|K| =m$ and $|I|$ even, $\alpha$ is a partition of $I$ into pairs and $\beta \in \{1, \ldots, \rho\}^J$. We set $\tau_{I, \alpha}$ to be the product of $\tau_{i, j}$'s corresponding to pairs in $\alpha$, $l_{J, \beta} = \prod_{j \in J} l^{\beta(j)}_j$ and $o_K = \prod_{k \in K}o_k$ (again, we set $\tau_{\emptyset, \alpha} = l_{\emptyset, \beta} = o_\emptyset = 1$). There is a bijection between $\Mon^m(n)$ and $\Mon^{2n - m}(n)$ given by
\begin{equation*}
\tau_{I, \alpha} \cdot l_{J, \beta} \cdot o_K \tofrom \tau_{I, \alpha} \cdot l_{J, \beta} \cdot o_{(I \cup J \cup K)^\complement}.
\end{equation*}

For $0 \leq m \leq 2n$, consider the pairing
\begin{equation} \label{eq:pair}
\bQ \cdot \Mon^m(n) \times \bQ \cdot \Mon^{2n - m}(n) \to \bQ \cdot \Mon^{2n}(n) = \bQ \cdot \bigg\{\prod_{i = 1}^n o_i\bigg\},
\end{equation}
which is given by the same recipes \eqref{eq:bv1'}, \eqref{eq:bv2'} and \eqref{eq:bv3'}. We have the following observation.

\subsection{Lemma} \label{lem:mehdi} --- \textit{Let $\tau_{I, \alpha} \cdot l_{J, \beta} \cdot o_K$ and $\tau_{I', \alpha'} \cdot l_{J', \beta'} \cdot o_{K'}$ be elements in $\Mon^m(n)$ and $\Mon^{2n - m}(n)$ respectively. Then the pairing of the two can be non-zero only if $I' = I, J' = J, K' = (I \cup J \cup K)^\complement$ and $\beta' = \beta$.}

\begin{proof}
Suppose the pairing is non-zero. We have
\begin{equation*}
2n = |I| + |J| + 2|K| + |I'| + |J'| + 2|K'| \leq 2|I \cup I'| + 2|J \cup J'| + 2|K| + 2|K'| \leq 2n.
\end{equation*}
Here the first inequality is obvious and the second follows from the fact that $I \cup I', J \cup J', K$ and $K'$ are pairwise disjoint, which in turn follows from \eqref{eq:bv1'}, \eqref{eq:bv2'} and \eqref{eq:bv3'}. Therefore the two inequalities are both equalities, which implies $I' = I$, $J' = J$ and $K' = (I \cup J \cup K)^\complement$. Further, the assumption that $\{l^s\}$ is an orthogonal basis implies $\beta' = \beta$.
\end{proof}

\subsection{\texorpdfstring{\!\!\!}{}} Lemma~\ref{lem:mehdi} shows that after a suitable ordering of the bases, the pairing matrix of \eqref{eq:pair} is block diagonal. Moreover, the diagonal blocks correspond to the pairing of elements in $\Mon^d(d)$ for some $d \leq \min\{m, 2n - m\}$ that consist solely of $\tau_{i, j}$'s. We denote by $\Mon_\tau^d(d) \subset \Mon^d(d)$ the subset of all such elements, \ie monomials of the form $\tau_\alpha = \tau_{\{1, \ldots, d\}, \alpha}$ where $\alpha$ is a partition of $\{1, \ldots, d\}$ into pairs ($d$ even).

We are left to consider the pairing
\begin{equation} \label{eq:pair2}
\bQ \cdot \Mon_\tau^d(d) \times \bQ \cdot \Mon_\tau^d(d) \to \bQ \cdot \Mon^{2d}(d) = \bQ \cdot \bigg\{\prod_{i = 1}^d o_i\bigg\}.
\end{equation}
It turns out that the matrix of \eqref{eq:pair2} has been studied in details by Hanlon and Wales \cite{HW89}. It is denoted by $T_r(x)$ with $r = d/2$ and $x = b_{\tr}$. Here we only cite (and translate) what is relevant to our problem, namely Theorem~3.1 in \loccit.

\subsection{Proposition} \label{prop:hw} --- \textit{The symmetric group $\fS_d$ acts on $\bQ \cdot \Mon_\tau^d(d)$ by permutations. We have
\begin{equation*}
\bQ \cdot \Mon_\tau^d(d) \simeq \bigoplus_{\lambda \in \Lambda_d} V_\lambda,
\end{equation*}
where $\Lambda_d$ is the set of partitions of $\{1, \ldots, d\}$ whose parts are all even, and $V_\lambda$ is the irreducible representation of $\fS_d$ associated to $\lambda$. All $V_\lambda$'s are eigenspaces of the matrix $T_{d/2}(b_{\tr})$. The eigenvalue is $0$ if and only if $\lambda$ contains at least $b_{\tr} + 1$ parts.}

\subsection{End of proof} --- First recall the definition of $V_\lambda$ via Specht modules (see \cite{FH91}, Problem~4.47). Define a tabloid $\{T\}$ to be an equivalence class of Young tableaux $T$ associated to $\lambda$, two being equivalent if the rows are the same up to order. The group $\fS_d$ acts by permutations on the set of such tabloids, denoted by $\Tab(\lambda)$. For each tableau $T$, define
\begin{equation*}
E_T = \sum_{g \in Q_T}\sgn(g)\big\{g(T)\big\} \in \bQ \cdot \Tab(\lambda),
\end{equation*}
where $Q_T \subset \fS_d$ is the column stabilizer of $T$. Then $V_\lambda$ is the $\bQ$-span of all $E_T$'s in $\bQ \cdot \Tab(\lambda)$. A basis of $V_\lambda$ is given by the $E_T$'s with standard tableaux $T$.

Now we locate $V_\lambda$ inside $\bQ \cdot \Mon_\tau^d(d)$. Denote by $T_i \subset \{1, \ldots, d\}$ the $i$-th row of a tableau $T$. For each $T_i$ consider the sum $\sum_{\alpha_i} \tau_{T_i, \alpha_i}$, where $\alpha_i$ runs through all partitions of $T_i$ into pairs. Altogether we assign $T$ the product
\begin{equation*}
\phi(T) = \prod_{i}\Big(\sum_{\alpha_i} \tau_{T_i, \alpha_i}\Big) \in \bQ \cdot \Mon_\tau^d(d).
\end{equation*}
It is easy to see that $\phi$ passes to the tabloids and that the resulting map $\phi \colon \bQ \cdot \Tab(\lambda) \to \bQ \cdot \Mon_\tau^d(d)$ is $\fS_d$-equivariant. By restriction we get a morphism of $\fS_d$-modules $\phi|_{V_\lambda} \colon V_\lambda \to \bQ \cdot \Mon_\tau^d(d)$, which is injective since it is non-zero and $V_{\lambda}$ is irreducible. We may identify $V_\lambda$ with its image.

By Proposition~\ref{prop:hw}, we know that $V_\lambda$ lies in the kernel of \eqref{eq:pair2} if and only if $\lambda$ contains at least $b_{\tr} + 1$ parts. The first occurrence is when $d = 2(b_{\tr} + 1)$ and $\lambda = (2, \ldots, 2)$. Take for example the following standard tableau $T$, with $Q_T \simeq \fS_{b_{\tr} + 1} \times \fS_{b_{\tr} + 1}$.
\begin{center} \medskip
\ytableausetup{mathmode,boxsize=2.1em}
\begin{ytableau}
\scriptstyle 1 & \scriptscriptstyle b_{\tr} + 2 \\
\scriptstyle 2 & \scriptscriptstyle b_{\tr} + 3 \\
\none[\vdots] & \none[\vdots] \\
\scriptscriptstyle b_{\tr} + 1 & \scriptscriptstyle 2b_{\tr} + 2
\end{ytableau} \medskip
\end{center}
A direct calculation gives
\begin{equation} \label{eq:motrel2}
\phi(E_T) = (b_{\tr} + 1)! \sum_{g \in \fS_{b_{\tr} + 1}} \sgn(g) \prod_{i = 1}^{b_{\tr} + 1} \tau_{i, b_{\tr} + 1 + g(i)},
\end{equation}
which is exactly $(b_{\tr} + 1)!$ times the left-hand side of \eqref{eq:motrel}. The other $\phi(E_T)$'s are given by permuting the indices on the right-hand side of \eqref{eq:motrel2}. 

The situation is similar as long as $\lambda$ contains at least $b_{\tr} + 1$ parts. We draw a standard tableau~$T$ of such a $\lambda$, with the length of the first and second columns $e \geq b_{\tr} + 1$ (lengths of other columns do not matter). Here $Q_\lambda \simeq \fS_e \times \fS_e \times \cdots$.
\begin{center} \medskip
\ytableausetup{mathmode,boxsize=2.1em}
\begin{ytableau}
\scriptstyle 1 & \scriptstyle e + 1 & \none[\cdots] & \scriptstyle d - 3 & \scriptstyle d - 1\\
\scriptstyle 2 & \scriptstyle e + 2 & \none[\cdots] & \scriptstyle d - 2 & \scriptstyle d\\
\none[\vdots] & \none[\vdots] \\
\scriptstyle e & \scriptstyle 2e
\end{ytableau} \medskip
\end{center}
By writing $\fS_e$ as the union of cosets $\{g \cdot \fS_{b_{\tr} + 1}\}$, it is not difficult to see that $\phi(E_T)$ is generated by various pull-backs and permutations of the right-hand side of \eqref{eq:motrel2}. The same holds for the other $\phi(E_T)$'s again by permutations.

We conclude that the kernel of \eqref{eq:pair2} is entirely generated by the right-hand side of \eqref{eq:motrel2} and its permutations. Then by the assumption \eqref{eq:motrel} those classes vanish in $R(S^{2(b_{\tr} + 1)})$. It follows that the kernel of \eqref{eq:pair2} vanishes in $R(S^d)$, and that the pairing between $R^m(S^n)$ and $R^{2n - m}(S^n)$ is perfect. \qed

\subsection{Final remarks} --- (i) \ A $1$-dimensional analogue of our result concerns a hyperelliptic curve~$C$. There the ring $R(C^n)$ contains (pull-backs of) the canonical class of $C$ and the diagonal class on $C \times C$. It is proven that all relations in $R(C^n)$ are generated by the vanishing of the Faber-Pandharipande cycle in $R^2(C^2)$, the vanishing of the Gross-Schoen cycle in $R^2(C^3)$, and one relation corresponding to the vanishing of the motive $\ssS^{2g + 2}h^1(C)$ (unconditional). This is the work of Tavakol (\cite{Tav11}; see also \cite{Tav14}), which inspired the present note. It would also be interesting to see if there are higher-dimensional analogues.

(ii) \ As is remarked in \cite{Voi14}, Section~5.1, it might be the case that Conjecture~\hyperref[conj:bv]{2} for $S^{[n]}$ also implies Conjecture~\hyperref[conj:v]{1}, \ie the two are equivalent. Our result shows that it suffices to deduce the relation \eqref{eq:motrel} from Conjecture~\hyperref[conj:bv]{2}. By the work of de Cataldo and Migliorini on the decomposition of $\CH(S^{[n]})$ (\cite{CM02}, Theorem~5.4.1), one can express the left-hand side of \eqref{eq:motrel} as a homologically trivial class in $\CH(S^{[n]})$ for some very large $n$ ($n \geq 1 + 2 + \cdots + (b_{\tr} + 1)$). It remains to see if this class is generated by divisor classes and Chern classes of the tangent bundle. The computation feels like a ``reverse engineering'' of \cite{Voi08}, proof of Proposition~2.6.

(iii) \ Further, observe that the relations \eqref{eq:bv1'}, \eqref{eq:bv2'}, \eqref{eq:bv3'} and \eqref{eq:motrel} hold in cohomology for any smooth projective surface of Albanese dimension~$0$. The same argument then shows that for such a surface~$S$, all relations in $\cl\big(R(S^n)\big)$ are generated by \eqref{eq:bv1'}, \eqref{eq:bv2'}, \eqref{eq:bv3'} and \eqref{eq:motrel}. We refer to \cite{Voi14b} and \cite{Via14} for recent applications of this result.

\secskip

\end{document}